\documentclass[11pt]{article}
\usepackage{bbm}
\usepackage{mathrsfs}
\usepackage{amscd}
\usepackage{amsmath,amsfonts,amssymb,amscd}
\usepackage{indentfirst,graphics,epsfig,psfrag}
\input{epsf}
\usepackage{ifpdf}
\usepackage{enumerate}
\usepackage{appendix}
\usepackage{enumerate}

\usepackage{lineno}
\makeatletter \@addtoreset{figure}{section} \makeatother
\makeatletter
\long\def\@makecaption#1#2{%
   \vskip 10\p@
   \setbox\@tempboxa\hbox{{#1}\ \ #2}%
   \ifdim \wd\@tempboxa >\hsize

       {#1}\ \ #2\par
   \else
       \hbox to\hsize{\hfil\box\@tempboxa\hfil}%
   \fi}
\makeatother

\newtheorem{thm}{Theorem}[section]
\newtheorem{cor}{Corollary}[section]
\newtheorem{lem}{Lemma}[section]

\newtheorem{obs}{Observation}[section]
\newtheorem{pro}{Proposition}[section]

\newcommand{\qed}{{\hfill\rule{3pt}{7pt}}}

\def\qed{\hfill \rule{4pt}{7pt}}

\begin{document}
\title{\textbf{Steiner diameter, maximum degree and size\\ of a graph}
\footnote{Supported by the National Science Foundation of China
(Nos. 11601254, 11551001, 11161037, and 11461054) and the Science Found of Qinghai
Province (Nos. 2016-ZJ-948Q, and 2014-ZJ-907).}}
\author{
\small Yaping Mao$^{1,2}$, \ \ Zhao Wang$^3$\\[0.2cm]
\small $^1$Department of Mathematics, Qinghai Normal\\
\small University, Xining, Qinghai 810008, China\\[0.2cm]
\small $^2$Center for Mathematics and Interdisciplinary Sciences\\
\small  of Qinghai Province, Xining, Qinghai 810008, China\\[0.2cm]
\small $^3$School of Mathematical Sciences, Beijing Normal\\
\small University, Beijing 100875, China}
\date{}
\maketitle
\begin{abstract}
The Steiner diameter $sdiam_k(G)$ of a graph $G$, introduced by
Chartrand, Oellermann, Tian and Zou in 1989, is a natural
generalization of the concept of classical diameter. When $k=2$,
$sdiam_2(G)=diam(G)$ is the classical diameter. The problem of
determining the minimum size of a graph of order $n$ whose diameter
is at most $d$ and whose maximum is $\ell$ was first introduced by
Erd\"{o}s and R\'{e}nyi. Recently, Mao considered the problem of
determining the minimum size of a graph of order $n$ whose Steiner
$k$-diameter is at most $d$ and whose maximum is at most $\ell$,
where $3\leq k\leq n$, and studied this new problem when $k=3$. In
this paper, we
investigate the problem when $n-3\leq k\leq n$.\\[2mm]
\end{abstract}

\section{Introduction}

All graphs in this paper are undirected, finite and simple. We refer
to \cite{Bondy} for graph theoretical notation and terminology not
described here. For a graph $G$, let $V(G)$, $E(G)$, $e(G)$,
$\delta(G)$, and $\overline{G}$ denote the set of vertices, the set
of edges, the size, minimum degree, and the complement of $G$,
respectively. The \emph{connectivity} $\kappa(G)$ is defined as the
order of a minimum vertex subset $S$ of $V(G)$ such that $G-S$ is
disconnected or has only one vertex. In this paper, we let $K_{n}$,
$P_n$, $K_{1,n-1}$ and $C_n$ be the complete graph of order $n$, the
path of order $n$, the star of order $n$, and the cycle of order
$n$, respectively. For any subset $X$ of $V(G)$, let $G[X]$ denote
the subgraph induced by $X$; similarly, for any subset $F$ of
$E(G)$, let $G[F]$ denote the subgraph induced by $F$. We use
$G\setminus X$ to denote the subgraph of $G$ obtained by removing
all the vertices of $X$ together with the edges incident with them
from $G$; similarly, we use $G\setminus F$ to denote the subgraph of
$G$ obtained by removing all the edges of $F$ from $G$. For two
subsets $X$ and $Y$ of $V(G)$ we denote by $E_G[X,Y]$ the set of
edges of $G$ with one end in $X$ and the other end in $Y$. The {\it
union}\index{union} $G\cup H$ of two graphs $G$ and $H$ is the graph
with vertex set $V(G)\cup V(H)$ and edge set $E(G)\cup E(H)$. If $G$
is the disjoint union of $k$ copies of a graph $H$, we simply write
$G=kH$. The {\it join}\index{join} $G\vee H$ of two disjoint graphs
$G$ and $H$ is the graph with vertex set $V(G)\cup V(H)$ and edge
set $E(G)\cup E(H)\cup \{uv\,|\, u\in V(G), v\in V(H)\}$. We divide
our introduction into the following four subsections to state the
motivations and our results of this paper.

\subsection{Distance and its generalizations}

Distance is one of the most basic concepts of graph-theoretic
subjects. If $G$ is a connected graph and $u,v\in V(G)$, then the
\emph{distance} $d_G(u,v)$ between $u$ and $v$ is the length of a
shortest path connecting $u$ and $v$. If $v$ is a vertex of a
connected graph $G$, then the \emph{eccentricity} $e(v)$ of $v$ is
defined by $e(v)=\max\{d_G(u,v)\,|\,u\in V(G)\}$. Furthermore, the
\emph{radius} $rad(G)$ and \emph{diameter} $diam(G)$ of $G$ are
defined by $rad(G)=\min\{e(v)\,|\,v\in V(G)\}$ and $diam(G)=\max
\{e(v)\,|\,v\in V(G)\}$. These last two concepts are related by the
inequalities $rad(G)\leq diam(G) \leq 2 rad(G)$. The \emph{center}
$C(G)$ of a connected graph $G$ is the subgraph induced by the
vertices $u$ of $G$ with $e(u)=rad(G)$. Goddard and Oellermann gave
a survey on this subject, see \cite{Goddard}.

The distance between two vertices $u$ and $v$ in a connected graph
$G$ also equals the minimum size of a connected subgraph of $G$
containing both $u$ and $v$. This observation suggests a
generalization of distance. The Steiner distance of a graph,
introduced by Chartrand, Oellermann, Tian and Zou in 1989, is a
natural generalization of the concept of classical graph distance.
For a graph $G(V,E)$ and a set $S\subseteq V(G)$ of at least two
vertices, \emph{an $S$-Steiner tree} or \emph{a Steiner tree
connecting $S$} (or simply, \emph{an $S$-tree}) is a subgraph
$T(V',E')$ of $G$ that is a tree with $S\subseteq V'$. Let $G$ be a
connected graph of order at least $2$ and let $S$ be a nonempty set
of vertices of $G$. Then the \emph{Steiner distance} $d_G(S)$ among
the vertices of $S$ (or simply the distance of $S$) is the minimum
size among all connected subgraphs whose vertex sets contain $S$.
Note that if $H$ is a connected subgraph of $G$ such that
$S\subseteq V(H)$ and $|E(H)|=d_G(S)$, then $H$ is a tree. Observe
that $d_G(S)=\min\{e(T)\,|\,S\subseteq V(T)\}$, where $T$ is subtree
of $G$. Furthermore, if $S=\{u,v\}$, then $d_G(S)=d(u,v)$ is the
classical distance between $u$ and $v$. Set $d_G(S)=\infty$ when
there is no $S$-Steiner tree in $G$.
\begin{obs}\label{obs1-1}
Let $G$ be a graph of order $n$ and $k$ be an integer with $2\leq
k\leq n$. If $S\subseteq V(G)$ and $|S|=k$, then $d_G(S)\geq k-1$.
\end{obs}

Let $n$ and $k$ be two integers with $2\leq k\leq n$. The
\emph{Steiner $k$-eccentricity $e_k(v)$} of a vertex $v$ of $G$ is
defined by $e_k(v)=\max \{d(S)\,|\,S\subseteq V(G), |S|=k,~and~v\in
S \}$. The \emph{Steiner $k$-radius} of $G$ is $srad_k(G)=\min \{
e_k(v)\,|\,v\in V(G)\}$, while the \emph{Steiner $k$-diameter} of
$G$ is $sdiam_k(G)=\max \{e_k(v)\,|\,v\in V(G)\}$. Note for every
connected graph $G$ that $e_2(v)=e(v)$ for all vertices $v$ of $G$
and that $srad_2(G)=rad(G)$ and $sdiam_2(G)=diam(G)$.

The following Table 1 shows how the generalization proceeds. {\small
\begin{center}
\begin{tabular}{|c|c|c|}
\hline  & Classical distance parameters& Steiner distance parameters\\[0.1cm]
\cline{1-3}
Vertex subset & $S=\{x,y\}\subseteq V(G) \ (|S|=2)$ & $S\subseteq V(G) \ (|S|=k\geq 2)$\\[0.1cm]
\cline{1-3} (Steiner-) distance& $\left\{
\begin{array}{ll}
d_G(x,y)=\min_{\{x,y\}\subseteq V(H)}\{e(H)\}\\[0.1cm]
H~is~a~subgraph~of~G\\[0.1cm]
\end{array}
\right.$ & $\left\{
\begin{array}{ll}
d_G(S)=\min_{S\subseteq V(H)}\{e(H)\}\\[0.1cm]
H~is~a~subgraph~of~G\\[0.1cm]
\end{array}
\right.$\\[0.1cm]
\cline{1-3} (Steiner) eccentricity & $e(v)=\max \{d_G(x,y)\,|\,~v\in
\{x,y\}\}$ & $\left\{
\begin{array}{ll}
e_k(v)=\max \{d_G(S)\,|\,~v\in
S \}\\[0.1cm]
S\subseteq V(G), |S|=k\\[0.1cm]
\end{array}
\right.$\\[0.1cm]
\cline{1-3} (Steiner) radius & $srad(G)=\min \{ e(v)\,|\,v\in
V(G)\}$ & $srad_k(G)=\min \{
e_k(v)\,|\,v\in V(G)\}$\\[0.1cm]
\cline{1-3}
(Steiner) diameter & $diam(G)=\max \{e(v)\,|\,v\in V(G)\}$ & $sdiam_k(G)=\max \{e_k(v)\,|\,v\in V(G)\}$\\[0.1cm]
\cline{1-3}
\end{tabular}
\end{center}
\begin{center}
{Table 1. Classical~distance parameters and Steiner distance
parameters}
\end{center}
}

\begin{obs}\label{obs1-2}
Let $k,n$ be two integers with $2\leq k\leq n$.

$(1)$ If $H$ is a spanning subgraph of $G$, then $sdiam_k(G)\leq
sdiam_k(H)$.

$(2)$ For a connected graph $G$, $sdiam_k(G)\leq sdiam_{k+1}(G)$.
\end{obs}

In \cite{ChartrandOZ}, Chartrand, Okamoto, Zhang obtained the
following upper and lower bounds of $sdiam_k(G)$.
\begin{thm}{\upshape\cite{ChartrandOZ}}\label{th1-1}
Let $k,n$ be two integers with $2\leq k\leq n$, and let $G$ be a
connected graph of order $n$. Then $k-1\leq sdiam_k(G)\leq n-1$.
Moreover, the upper and lower bounds are sharp.
\end{thm}

In \cite{DankelmannSO2}, Dankelmann, Swart and Oellermann obtained a
bound on $sdiam_k(G)$ for a graph $G$ in terms of the order of $G$
and the minimum degree $\delta$ of $G$, that is, $sdiam_k(G)\leq
\frac{3n}{\delta+1}+3k$. Later, Ali, Dankelmann, Mukwembi
\cite{AliDM} improved the bound of $sdiam_k(G)$ and showed that
$sdiam_k(G)\leq \frac{3n}{\delta+1}+2k-5$ for all connected graphs
$G$. Moreover, they constructed graphs to show that the bounds are
asymptotically best possible. In \cite{Mao}, Mao obtained the
Nordhaus-Gaddum-type results for the parameter $sdiam_k(G)$.

As a generalization of the center of a graph, the \emph{Steiner
$k$-center} $C_k(G)\ (k\geq 2)$ of a connected graph $G$ is the
subgraph induced by the vertices $v$ of $G$ with $e_k(v)=srad_k(G)$.
Oellermann and Tian \cite{OellermannT} showed that every graph is
the $k$-center of some graph. In particular, they showed that the
$k$-center of a tree is a tree and those trees that are $k$-centers
of trees are characterized. The \emph{Steiner $k$-median} of $G$ is
the subgraph of $G$ induced by the vertices of $G$ of minimum
Steiner $k$-distance. For Steiner centers and Steiner medians, we
refer to \cite{Oellermann, Oellermann2, OellermannT}.

The \emph{average Steiner distance} $\mu_k(G)$ of a graph $G$,
introduced by Dankelmann, Oellermann and Swart in
\cite{DankelmannOS}, is defined as the average of the Steiner
distances of all $k$-subsets of $V(G)$, i.e.
$$
\mu_k(G)={n\choose k}^{-1}\sum_{S\subseteq V(G), \ |S|=k}d_G(S).
$$
For more details on average Steiner distance, we refer to
\cite{DankelmannOS, DankelmannSO}.

Let $G$ be a $k$-connected graph and $u$, $v$ be any pair of
vertices of $G$. Let $P_k(u,v)$ be a family of $k$ inner
vertex-disjoint paths between $u$ and $v$, i.e.,
$P_k(u,v)=\{P_1,P_2,\cdots,P_k\}$, where $p_1\leq p_2\leq \cdots
\leq p_k$ and $p_i$ denotes the number of edges of path $P_i$. The
\emph{$k$-distance} $d_k(u,v)$ between vertices $u$ and $v$ is the
minimum $p_k$ among all $P_k(u,v)$ and the \emph{$k$-diameter}
$d_k(G)$ of $G$ is defined as the maximum $k$-distance $d_k(u,v)$
over all pairs $u,v$ of vertices of $G$. The concept of $k$-diameter
emerges rather naturally when one looks at the performance of
routing algorithms. Its applications to network routing in
distributed and parallel processing are studied and discussed by
various authors including Chung \cite{Chung}, Du, Lyuu and Hsu
\cite{Du}, Hsu \cite{Hsu, Hsu2}, Meyer and Pradhan \cite{Meyer}.

\subsection{Application background of Steiner distance parameters}

The Steiner tree problem in networks, and particularly in graphs,
was formulated in 1971-by Hakimi (see \cite{Hakimi}) and Levi (see
\cite{Levi}). In the case of an unweighted, undirected graph, this
problem consists of finding, for a subset of vertices $S$, a
minimal-size connected subgraph that contains the vertices in $S$.
The computational side of this problem has been widely studied, and
it is known that it is an NP-hard problem for general graphs (see
\cite{HwangRW}). The determination of a Steiner tree in a graph is a
discrete analogue of the well-known geometric Steiner problem: In a
Euclidean space (usually a Euclidean plane) find the shortest
possible network of line segments interconnecting a set of given
points. Steiner trees have application to multiprocessor computer
networks. For example, it may be desired to connect a certain set of
processors with a subnetwork that uses the least number of
communication links. A Steiner tree for the vertices, corresponding
to the processors that need to be connected, corresponds to such a
desired subnetwork.

The \emph{Wiener index} $W(G)$ of the graph $G$ is defined as
$W(G)=\sum_{\{u,v\} \subseteq V(G)} d_G(u,v)$. Details on this
oldest distance--based topological index can be found in numerous
surveys, e.g., in \cite{Dobrynin, Rouv1, Rouv2, Xu}. Li et al.
\cite{LMG} put forward a Steiner--distance--based generalization of
the Wiener index concept. According to \cite{LMG}, the {\it
$k$-center Steiner Wiener index\/} $SW_k(G)$ of the graph $G$ is
defined by
\begin{equation}                    \label{sw}
SW_k(G)=\sum_{\overset{S\subseteq V(G)}{|S|=k}} d(S)\,.
\end{equation}
For $k=2$, the above defined Steiner Wiener index coincides with the
ordinary Wiener index. It is usual to consider $SW_k$ for $2 \leq k
\leq n-1$, but the above definition would be applicable also in the
cases $k=1$ and $k=n$, implying $SW_1(G)=0$ and $SW_n(G)=n-1$. A
chemical application of $SW_k$ was recently reported in \cite{GFL}.
Gutman \cite{GutmanSDD} offered an analogous generalization of the
concept of degree distance. Later, Furtula, Gutman, and Katani\'{c}
\cite{FurtulaGK} introduced the concept of Steiner Harary index and
gave its chemical applications. Recently, Mao and Das \cite{MaoDas}
introduced the concept of Steiner Gutman index and obtained some
bounds for it. For more details on Steiner distance indices, we
refer to \cite{FurtulaGK, GFL, GutmanSDD, LMG, LMG2, MaoDas, MWG,
MWGK, MWGL}.

\subsection{Classical extremal problem and our generalization}

What is the minimal size of a graph of order $n$ and diameter $d$ ?
What is the maximal size of a graph of order $n$ and diameter $d$ ?
It is not surprising that these questions can be answered without
the slightest effort (see \cite{Bollobas}) just as the similar
questions concerning the connectivity or the chromatic number of a
graph. The class of maximal graphs of order $n$ and diameter $d$ is
easy to describe and reduce every question concerning maximal graphs
to a not necessarily easy question about binomial coefficient, as in
\cite{HO1, HS5, O6, W11}. Therefore, the authors study the minimal
size of a graph of order $n$ and under various additional
conditions.

Erd\"{o}s and R\'{e}nyi \cite{ER4} introduced the following problem.
Let $d,\ell$ and $n$ be natural numbers, $d<n$ and $\ell<n$. Denote
by $\mathscr{H}(n,\ell,d)$ the set of all graphs of order $n$ with
maximum degree $\ell$ and diameter at most $d$. Put
$$
e(n,\ell,d)=\min\{e(G):G\in \mathscr{H}(n,\ell,d)\}.
$$
If $\mathscr{H}(n,\ell,d)$ is empty, then, following the usual
convention, we shall write $e(n,\ell,d)=\infty$. For more details on
this problem, we refer to \cite{Bollobas, B26, ER4, ERS1}.

Mao \cite{Mao2} considered the generalization of the above problem.
Let $d,\ell$ and $n$ be natural numbers, $d<n$ and $\ell<n$. Denote
by $\mathscr{H}_k(n,\ell,d)$ the set of all graphs of order $n$ with
maximum degree $\ell$ and $sdiam_k(G)\leq d$. Put
$$
e_k(n,\ell,d)=\min\{e(G):G\in \mathscr{H}_k(n,\ell,d)\}.
$$
If $\mathscr{H}_k(n,\ell,d)$ is empty, then, following the usual
convention, we shall write $e_k(n,\ell,d)=\infty$. From Theorem
\ref{th1-1}, we have $k-1\leq d\leq n-1$.

In \cite{Mao2}, Mao focused their attention on the case $k=3$, and
studied the exact value of $e_3(n,\ell,d)$ for $d=n-1,n-2,n-3,2,3$.
In this paper, we investigate another extreme case when $n-3\leq
k\leq n-1$, and give the exact values or upper and lower bounds of
$e_k(n,\ell,d)$ for $n-3\leq k\leq n-1$. For general $k \ (3\leq
k\leq n-1)$, $d \ (k-1\leq d\leq n-1)$ and $\ell \ (2\leq \ell\leq
n-1)$, we give upper and lower bounds of $e_k(n,\ell,d)$.

\section{The case $k=n,n-1$}

In the sequel, let $K_{s,t}$, $K_{n}$, $C_{n}$ and $P_n$ denote the
complete bipartite graph of order $s+t$ with part sizes $s$ and $t$,
complete graph of order $n$, cycle of order $n$ and path of order
$n$, respectively.

The following observation is immediate.
\begin{obs}{\upshape \cite{Mao}} \label{obs2-1}
$(1)$ For a cycle $C_n$, $sdiam_k(C_n)=\left
\lfloor\frac{n(k-1)}{k}\right\rfloor$;

$(2)$ For a complete graph $K_n$, $sdiam_k(K_n)=k-1$.
\end{obs}

The following result is easily proved in \cite{Mao2}.
\begin{lem}{\upshape \cite{Mao2}}\label{lem2-1}
For $2\leq \ell \leq n-1$ and $3\leq k \leq n$,
$$
e_k(n,\ell,n-1)=n-1.
$$
\end{lem}

For $k=n$, we know that $sdiam_n(G)=n-1$ for a connected graph $G$,
and hence $d=n-1$. From Lemma \ref{lem2-1}, the following result is
immediate.
\begin{pro}\label{pro2-1}
For $2\leq \ell \leq n-1$, $e_n(n,\ell,n-1)=n-1$.
\end{pro}

From now on, we assume that $k\leq n-1$.
\begin{pro}\label{pro2-2}
For $3\leq k \leq n$, $e_k(n,n-1,k)=n-1$.
\end{pro}
\begin{pf}
Let $G=K_{1,n-1}$ be a star of order $n$. Clearly,
$\Delta(G)=\ell=n-1$. Since $sdiam_k(G)=k$ and $e(G)=n-1$, it
follows that $e_k(n,n-1,k)\leq n-1$. On the other hand, since we
only consider connected graphs, it follows that $e(G)\geq n-1$ for a
connected graph $G$ is of order $n$. So $e_k(n,n-1,k)=n-1$. \qed
\end{pf}

\vskip 0.5cm

For $k=n-1$, we have $n-2\leq sdiam_{n-1}(G)\leq n-1$ by Theorem
\ref{th1-1}. So we only need to consider the case $d=n-1$ or
$d=n-2$. Note that $2\leq \ell \leq n-1$.
\begin{thm}\label{th2-1}
$(1)$ For $2\leq \ell \leq n-1$, $e_{n-1}(n,\ell,n-1)=n-1$.

$(2)$ For $2\leq \ell \leq n-1$, $e_{n-1}(n,\ell,n-2)=n+\ell-2$.
\end{thm}
\begin{pf}
$(1)$ The result follows from Lemma \ref{lem2-1}.

$(2)$ For $\ell=2$, we let $C_n$ be the cycle of order $n$. From
Observation \ref{obs2-1}, we have $sdiam_{n-1}(G)=n-2$. Since
$e(G)=n$, it follows that $e_{n-1}(n,2,n-2)\leq n$. Let $G$ be graph
of order $n$ with $\Delta(G)=2$ and $sdiam_{n-1}(G)=n-2$. Since
$\Delta(G)=2$, it follows that $G=P_n$ or $G=C_n$. If $G=P_n$, then
$sdiam_{n-1}(G)=n-1$, a contradiction. Therefore, $G=C_n$ and hence
$e(G)\geq n$. So $e_{n-1}(n,2,n-2)=n$.

Suppose $3\leq \ell \leq n-1$. Let $G$ be a graph obtained from a
cycle $C_n=v_1v_2\cdots v_{n}$ by adding the edges $v_1v_{j} \
(3\leq j\leq \ell)$. From Observation \ref{obs1-1}, we have
$sdiam_{n-1}(G)\leq sdiam_{n-1}(C_n)\leq n-2$. Since
$\Delta(G)=\ell$ and $e(G)=n+\ell-2$, it follows that
$e_{n-1}(n,\ell,n-2)\leq n+\ell-2$. Conversely, let $G$ be a graph
such that $sdiam_{n-1}(G)=n-2$ and $\Delta(G)=\ell \ (3\leq \ell
\leq n-1)$. Then there exists a vertex $u$ in $G$ such that
$d_G(u)=\ell$. Choose $S=V(G)-u$. Since $sdiam_{n-1}(G)=n-2$, it
follows that there exists an $S$-Steiner tree in $G-u$, say $T$.
Then $e(G-u)\geq e(T)=n-2$ and hence $e(G)=e(G-u)+\ell\geq
n-2+\ell$, which implies that $e_{n-1}(n,\ell,n-2)\geq n+\ell-2$.

From the above arguments, we conclude that $e_{n-1}(n,\ell,n-2)=
n+\ell-2$ for $2\leq \ell \leq n-1$. \qed
\end{pf}

\section{The case $k=n-2$}

From Theorem \ref{th1-1}, we have $n-3\leq sdiam_{n-2}(G)\leq n-1$.
Mao et al. \cite{MaoMelekianCheng} characterized the graphs with
$sdiam_{n-2}(G)=d \ (n-3\leq d\leq n-1)$.
\begin{lem}{\upshape \cite{MaoMelekianCheng}}\label{lem3-1}
Let $G$ be a connected graph of order $n \ (n\geq 5)$. Then

$(1)$ $sdiam_{n-2}(G)=n-3$ if and only if $\kappa(G)\geq 3$.

$(2)$ $sdiam_{n-2}(G)=n-2$ if and only if $\kappa(G)=2$ or $G$
contains only one cut vertex.

$(3)$ $sdiam_{n-2}(G)=n-1$ if and only if there are at least two cut
vertices in $G$.
\end{lem}

For $d=n-2$, we have the following.
\begin{pro}\label{pro3-1}
For $2\leq \ell \leq n-1$ and $n\geq 5$,
$$
e_{n-2}(n,\ell,n-2)=\left\{
\begin{array}{ll}
n,&\mbox{{\rm if}~$2\leq \ell \leq n-2$;}\\
n-1,&\mbox{{\rm if}~$\ell=n-1$.}
\end{array}
\right.
$$
\end{pro}
\begin{pf}
For $\ell=n-1$, from Proposition \ref{pro2-2}, we have
$e_{n-2}(n,\ell,n-2)=n-1$. From now on, we suppose $2\leq \ell \leq
n-2$. Let $G$ be a graph obtained by a cycle $C_{n-\ell+2}$ and a
star $K_{1,\ell-2}$ by identifying a vertex of $C_{n-\ell+2}$ and
the center of $K_{1,\ell-2}$. Clearly, $\Delta(G)=\ell$ and there is
exactly one cut vertex in $G$. From $(2)$ of Lemma \ref{lem3-1}, we
have $sdiam_{n-2}(G)=n-2$, and hence $e_{n-2}(n,\ell,n-2)\leq n$. It
suffices to show that $e_{n-2}(n,\ell,n-2)\geq n$. Let $G$ be a
graph of order $n$ with $sdiam_{n-2}(G)\leq n-2$ and
$\Delta(G)=\ell$. If $G$ is a tree, then $G$ contains at least two
cut vertices, since $\Delta(G)=\ell\leq n-2$. From $(3)$ of Lemma
\ref{lem3-1}, $sdiam_{n-2}(G)=n-1$, a contradiction. So $G$ contains
at least one cycle, and hence $e(G)\geq n$. Therefore, we have
$e_{n-2}(n,\ell,n-2)=n$, as desired. \qed
\end{pf}

\vskip 0.4cm

Let $P_j^i$ be a path of order $j$, where $1\leq i\leq r+2$. We call
the graph $K_1\vee (K_1\cup P_j^i)$ as a {\it
$(u_i,v_i,P_j^i)$-Fan}; see Figure 1 $(a)$. For $1\leq i\leq r$, we
choose $(u_i,v_i,P_2^i)$-Fan, and choose
$(u_{r+1},v_{r+1},P_{\ell-1}^{r+1})$-Fan and
$(u_{r+2},v_{r+2},P_{s}^{r+2})$-Fan. Let $H_n$ be a graph obtained
from the above $(r+2)$ Fans by adding the edges in
\begin{eqnarray*}
&&\{w_1^{r+1}w_{s}^{r+2},w_1^{r+2}w_{1}^{1}\}\cup \{w_2^{i}w_{1}^{i+1}\,|\,1\leq i\leq r-1\}\cup \{w_2^{r}w_{\ell-1}^{r+1}\}\\[0.1cm]
&&\cup \{v_{i}v_{i+1}\,|\,1\leq i\leq r+1\}\cup \{v_{r+2}v_{1}\};
\end{eqnarray*}
see Figure 1 $(b)$, where $4r+\ell+s+3=n$, $2\leq s\leq 5$ and $1\leq r\leq \frac{n-\ell-4}{4}$.
\begin{figure}[!hbpt]
\begin{center}
\includegraphics[scale=0.6]{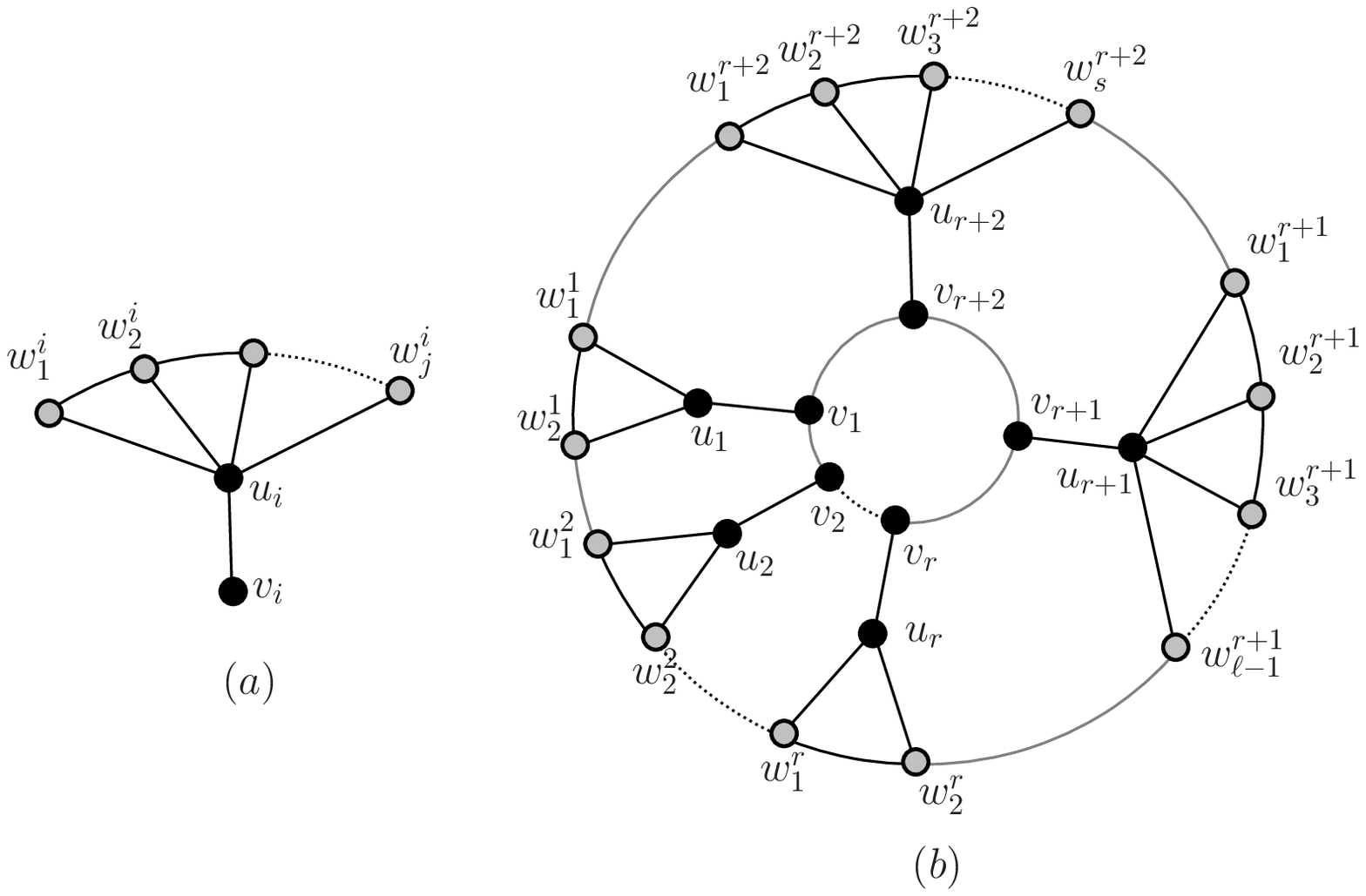}\\
Figure 1: Graphs for Proposition \ref{pro3-2}.
\end{center}\label{fig7}
\end{figure}

\begin{pro}\label{pro3-2}
For $6\leq \ell \leq n-9$,
$$
\frac{1}{2}(3n+\ell-3)\leq e_{n-2}(n,\ell,n-3)\leq \frac{1}{2}(3n+\ell+s-5),
$$
where $2\leq s\leq 5$. Furthermore, if $s=2$, then $e_{n-2}(n,\ell,n-3)=\frac{1}{2}(3n+\ell-3)$.
\end{pro}
\begin{pf}
Let $H_n$ be the graph constructed in Figure 1 $(b)$. Clearly, $H_n$
is $3$-connected. From Lemma \ref{lem3-1}, $sdiam_{n-2}(H_n)=n-3$.
Since $6\leq \ell\leq n-9$ and $2\leq s\leq 5$, it follows that
$\Delta(H_n)=\ell$. Then
\begin{eqnarray*}
e(G)&=&2s+2(\ell-1)+4r+(r-1)+r+2+3\\[0.1cm]
&=&6r+2s+2\ell+2\\[0.1cm]
&=&2\ell+(n-\ell-3)+2r+s+2\\[0.1cm]
&=&n+\ell+s-1+\frac{n-\ell-3-s}{2}\\[0.1cm]
&=&\frac{1}{2}\left[3n+\ell+s-5\right],
\end{eqnarray*}
and hence $e_{n-2}(n,\ell,n-3)\leq \frac{1}{2}(3n+\ell+s-5)$.

Conversely, we suppose that $G$ is a graph with $|V(G)|=n$, $\Delta(G)=\ell$, and $sdiam_{n-2}(G)=n-3$.
Then there is a vertex in $G$, say $u$, such that $d_G(u)=\ell$. Since $sdiam_{n-2}(G)=n-3$, it follows
from Lemma \ref{lem3-1} that $\kappa(G)\geq 3$, and hence $d_G(v)\geq 3$ for any $v\in V(G)\setminus u$.
Therefore, $e(G)\geq \frac{1}{2}(3n+\ell-3)$, and hence $e_{n-2}(n,\ell,n-3)\geq \frac{1}{2}(3n+\ell-3)$.\qed
\end{pf}

\vskip 0.3cm

For $n-8\leq \ell \leq n-1$, we have the following result.
\begin{lem}\label{lem3-2}
Let $k,\ell$ be two integers with $n-8\leq \ell \leq n-1$. Then

$(i)$ $e_{n-2}(n,n-1-i,n-3)=2n-2$ for $n\geq 5+i$ and $i=0,1$;

$(ii)$ $e_{n-2}(n,n-3-i,n-3)=2n-3$ for $n\geq 7+2i$ and $i=0,1$;

$(iii)$ $e_{n-2}(n,n-5-i,n-3)=2n-4$ for $n\geq 11+2i$ and $i=0,1$;

$(iv)$ $e_{n-2}(n,n-7-i,n-3)=2n-5$ for $n\geq 15+2i$ and $i=0,1$.
\end{lem}
\begin{pf}
For $(i)$, we first consider the case $i=0$. For $\ell=n-1$, let
$G_n^1$ be a wheel of order $n$. From Lemma \ref{lem3-1},
$sdiam_{n-2}(G_n^1)=n-3$ and $\Delta(G_n^1)=n-1$, and hence
$e_{n-2}(n,n-1,n-3)\leq 2n-2$. Conversely, we suppose that $G$ is a
graph of order $n$ such that $sdiam_{n-2}(G)=n-3$ and
$\Delta(G)=n-1$. Then there exists a vertex $u$ such that
$d_G(u)=n-1$. Since $\kappa(G)\geq 3$, it follows that
$\kappa(G-u)\geq 2$, and hence $e(G-u)\geq n-1$. Then $e(G)\geq
2n-2$, and hence $e_{n-2}(n,n-1,n-3)\geq 2n-2$. So
$e_{n-2}(n,n-1,n-3)=2n-2$.
\begin{figure}[!hbpt]
\begin{center}
\includegraphics[scale=0.7]{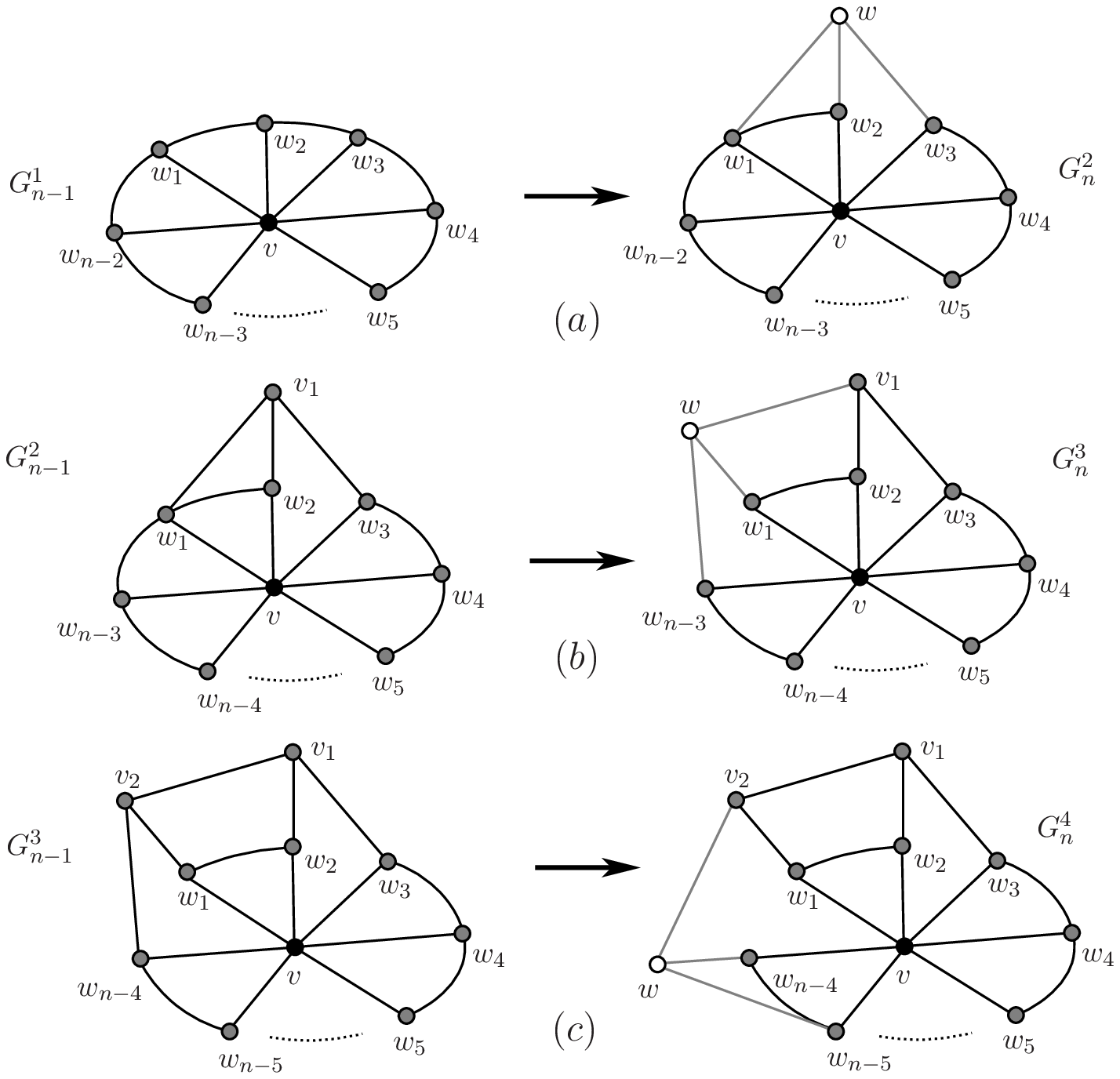}\\[0.5cm]
Figure 2: Graphs for $(i)$ of Lemma \ref{lem3-2}.
\end{center}\label{fig7}
\end{figure}
Next, we consider the case $i=1$. For $\ell=n-2$, let $G_n^2$ be a
graph of order $n$ obtained by $G_{n-1}^1$ by deleting the edge
$w_2w_3$, and then adding a new vertex $w$ and three edges
$ww_1,ww_2,ww_3$; see Figure 2 $(b)$. Since $\kappa(G_n^2)=3$, it
follows from Lemma \ref{lem3-1} that $sdiam_{n-2}(G_n^2)=n-3$. From
this together with $\Delta(G_n^2)=n-2$, we have
$e_{n-2}(n,n-2,n-3)\leq 2n-2$. Conversely, we suppose that $G$ is a
graph of order $n$ such that $sdiam_{n-2}(G)=n-3$ and
$\Delta(G)=n-2$. Then there exists a vertex $u$ such that
$d_G(u)=n-2$. Since $sdiam_{n-2}(G)=n-3$, it follows that
$\kappa(G)\geq 3$, and hence $\kappa(G-u)\geq 2$. Clearly,
$e(G-u)\geq n-1$. If $e(G-u)=n-1$, then $G-u$ is a cycle of order
$n-1$, say $G-u=v_1v_2\ldots v_{n-1}v_1$. Since $d_G(u)=n-2$, it
follows that there exists some vertex $v_i$ in $G-u$ such that
$uv_i\notin E(G)$, and hence $d_G(v_i)=2$, which contradicts to the
fact $\kappa(G)\geq 3$. Then $e(G-u)\geq n$, and hence $e(G)\geq
2n-2$, and hence $e_{n-2}(n,n-2,n-3)\geq 2n-2$. So, we have
$e_{n-2}(n,n-2,n-3)=2n-2$.
\begin{figure}[!hbpt]
\begin{center}
\includegraphics[scale=0.7]{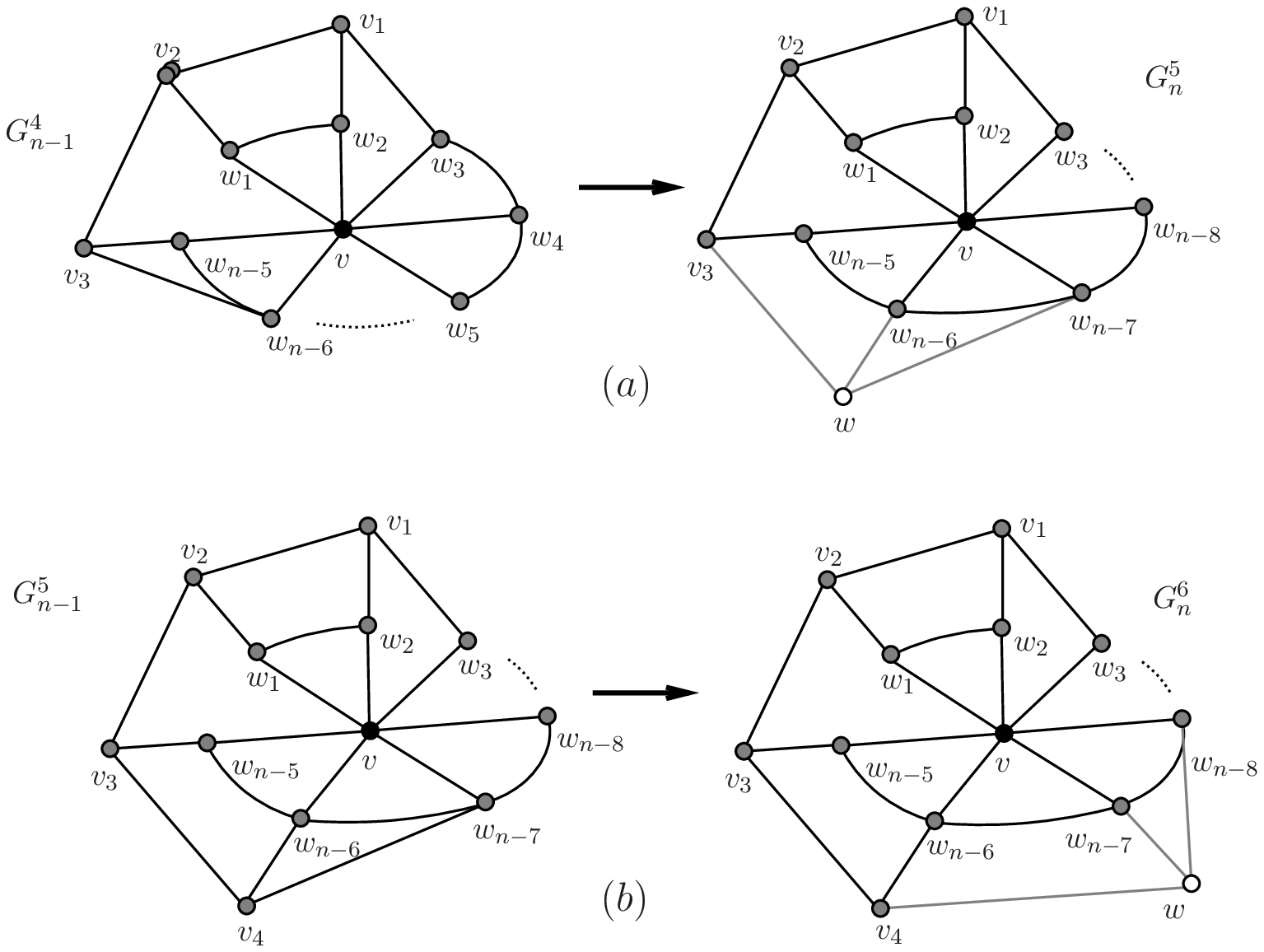}\\[0.5cm]
\includegraphics[scale=0.7]{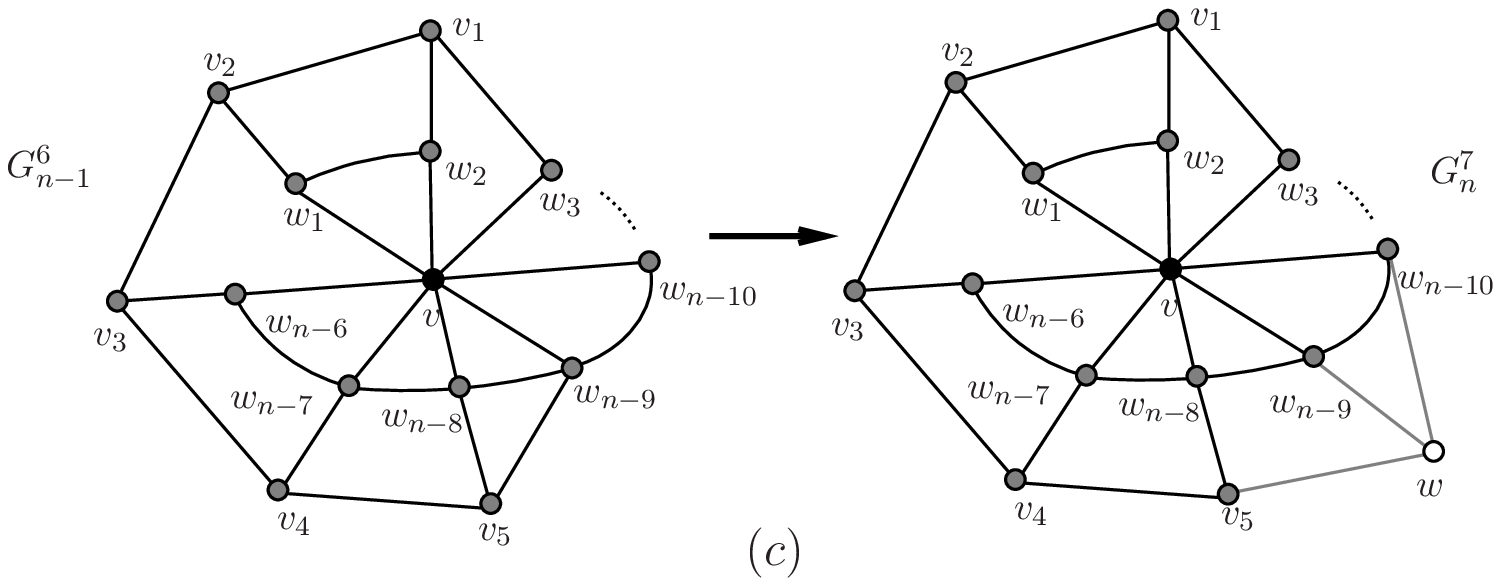}\\[0.5cm]
Figure 3: Graphs for $(iii)$ and $(iv)$ of Lemma \ref{lem3-2}.
\end{center}\label{fig7}
\end{figure}

For $(ii)$, we first consider the case $i=0$. For $\ell=n-3$, let
$G_n^3$ be a graph of order $n$ obtained by $G_{n-1}^2$ by deleting
the edge $v_1w_1,w_1w_{n-3}$, and then adding a new vertex $w$ and
three edges $wv_1,ww_1,ww_{n-3}$; see Figure 3 $(b)$. Since
$\kappa(G_n^3)=3$, it follows from Lemma \ref{lem3-1} that
$sdiam_{n-2}(G_n^3)=n-3$. Note that $\Delta(G_n^3)=n-3$. Therefore,
we have $e_{n-2}(n,n-3,n-3)\leq 2n-3$. Conversely, we suppose that
$G$ is a graph of order $n$ such that $sdiam_{n-2}(G)=n-3$ and
$\Delta(G)=n-3$. Then there exists a vertex $u$ such that
$d_G(u)=n-3$. Since $sdiam_{n-2}(G)=n-3$, it follows that
$\kappa(G)\geq 3$, and hence $\kappa(G-u)\geq 2$. Clearly,
$e(G-u)\geq n-1$. If $e(G-u)=n-1$, then $G-u$ is a cycle of order
$n-1$, say $G-u=v_1v_2\ldots v_{n-1}v_1$. Since $d_G(u)=n-3$, it
follows that there exists some vertex $v_i$ in $G-u$ such that
$uv_i\notin E(G)$, and hence $d_G(v_i)=2$, which contradicts to the
fact $\kappa(G)\geq 3$. Then $e(G-u)\geq n$, and hence $e(G)\geq
2n-3$, and hence $e_{n-2}(n,n-3,n-3)\geq 2n-3$. So, we have
$e_{n-2}(n,n-3,n-3)=2n-3$.

Next, we consider the case $i=1$. For $\ell=n-4$, let $G_n^4$ be a
graph of order $n$ obtained by $G_{n-1}^3$ by deleting the edge
$v_2w_{n-4}$, and then adding a new vertex $w$ and three edges
$wv_2,ww_{n-4},ww_{n-5}$; see Figure 3 $(d)$. Since
$\kappa(G_n^4)=3$, it follows from Lemma \ref{lem3-1} that
$sdiam_{n-2}(G_n^4)=n-3$. Note that $\Delta(G_n^3)=n-4$. Therefore,
we have $e_{n-2}(n,n-2,n-3)\leq 2n-3$. Conversely, we suppose that
$G$ is a graph of order $n$ such that $sdiam_{n-2}(G)=n-3$ and
$\Delta(G)=n-4$. Then there exists a vertex $u$ such that
$d_G(u)=n-4$. Since $sdiam_{n-2}(G)=n-3$, it follows that
$\kappa(G)\geq 3$, and hence $\kappa(G-u)\geq 2$. Clearly,
$e(G-u)\geq n-1$. If $e(G-u)=n-1$, then $G-u$ is a cycle of order
$n-1$, say $G-u=v_1v_2\ldots v_{n-1}v_1$. Since $d_G(u)=n-3$, it
follows that there exists some vertex $v_i$ in $G-u$ such that
$uv_i\notin E(G)$, and hence $d_G(v_i)=2$, which contradicts to the
fact $\kappa(G)\geq 3$. If $e(G-u)=n$, then $G-u$ is a graph
obtained from a cycle $C=v_1v_2\ldots v_{n-1}v_1$ by adding some
edge $v_pv_q \ (1\leq p\neq q\leq n-1)$. Since $d_G(u)=n-4$, it
follows that there exists some vertex $v_i \ (1\leq i\leq n-1, \
i\neq p, \ i\neq q)$ in $G-u$ such that $uv_i\notin E(G)$, and hence
$d_G(v_i)=2$, which contradicts to the fact $\kappa(G)\geq 3$. Then
$e(G-u)\geq n+1$, and hence $e(G)\geq 2n-3$, and hence
$e_{n-2}(n,n-4,n-3)\geq 2n-3$. So, we have
$e_{n-2}(n,n-4,n-3)=2n-3$.

For $(iii)$ and $(iv)$, we only give the graph construction
operation (see Figure 3), and omit the proof of them. \qed
\end{pf}

\vskip 0.3cm

From Propositions \ref{pro3-1} and \ref{pro3-2}, and Lemmas
\ref{lem3-2} and \ref{lem2-1}, we have the following theorem.
\begin{thm}\label{th3-2}
$(1)$ For $2\leq \ell \leq n-1$,
$$
e_{n-2}(n,\ell,n-1)=n-1.
$$

$(2)$ For $2\leq \ell \leq n-1$ and $n\geq 5$,
$$
e_{n-2}(n,\ell,n-2)=\left\{
\begin{array}{ll}
n,&\mbox{{\rm if}~$2\leq \ell \leq n-2$}\\
n-1,&\mbox{{\rm if}~$\ell=n-1$.}
\end{array}
\right.
$$

$(3)$ For $n-8\leq \ell \leq n-1$, $e_{n-2}(n,n-1-i,n-3)=2n-2$ for
$n\geq 5+i$ and $i=0,1$; $e_{n-2}(n,n-3-i,n-3)=2n-3$ for $n\geq
7+2i$ and $i=0,1$; $e_{n-2}(n,n-5-i,n-3)=2n-4$ for $n\geq 11+2i$ and
$i=0,1$; $e_{n-2}(n,n-7-i,n-3)=2n-5$ for $n\geq 15+2i$ and $i=0,1$.
For $6\leq \ell \leq n-9$,
$$
\frac{1}{2}(3n+\ell-3)\leq e_{n-2}(n,\ell,n-3)\leq
\frac{1}{2}(3n+\ell+s-5),
$$
where $2\leq s\leq 5$. Furthermore, if $s=2$, then
$e_{n-2}(n,\ell,n-3)=\frac{1}{2}(3n+\ell-3)$.
\end{thm}

\section{The case $k=n-3$}

Mao et al. \cite{Mao} derived the following results for Steiner
$(n-3)$-diameter.
\begin{lem}{\upshape \cite{Mao}}\label{lem4-1}
Let $G$ be a connected graph of order $n$. Then $sdiam_{n-3}(G)=n-4$
if and only if $\kappa(G)\geq 4$, and $sdiam_{n-3}(G)=n-1$ if and
only if $G$ contains at least $3$ cut vertices.
\end{lem}

Wang et al. \cite{WangMLY} obtained the structural properties of
graphs with $sdiam_k(G)=n-1$.
\begin{lem}{\upshape \cite{WangMLY}}\label{lem4-2}
Let $k,n$ be two integers with $3\leq k\leq n-1$. Let $G$ be a
connected graph of order $n$. Then $sdiam_k(G)=n-1$ if and only if
the number of non-cut vertices in $G$ is at most $k$.
\end{lem}

The following corollary is immediate from Lemma \ref{lem2-1}.
\begin{cor}\label{cor4-1}
For $2\leq \ell \leq n-1$, $e_{n-3}(n,\ell,n-1)=n-1$.
\end{cor}

Let $uv$ be an edge in $G$. A \emph{double-star} on $uv$ is a
maximal tree in $G$ which is the union of stars centered at $u$ or
$v$ such that each star contains the edge $uv$.
\begin{pro}\label{pro4-1}
For $2\leq \ell \leq n-1$ and $n\geq 4$,
$$
e_{n-3}(n,\ell,n-2)=\left\{
\begin{array}{ll}
n,&\mbox{{\rm if}~$2\leq \ell \leq \lfloor\frac{n}{2}\rfloor-1$};\\
&\mbox{{\rm ~~~or}~$\ell=\lfloor\frac{n}{2}\rfloor$~and~n~is~odd};\\
n-1,&\mbox{{\rm if}~$\lfloor\frac{n}{2}\rfloor+1\leq \ell \leq n-1$};\\
&\mbox{{\rm ~~~or}~$\ell=\lfloor\frac{n}{2}\rfloor$~and~n~is~even}.\\
\end{array}
\right.
$$
\end{pro}
\begin{pf}
For $2\leq \ell \leq \lfloor\frac{n}{2}\rfloor$, we let $G$ be a
graph of order $n$ obtained from a cycle $C_{n-\ell+2}$ and a star
$K_{1,\ell-2}$ by identifying the center and one vertex of
$C_{n-\ell+2}$. Clearly, $\Delta(G)=\ell$, and $G$ contains only one
cut vertex. From Lemma \ref{lem4-2}, we have $sdiam_{n-3}(G)\leq
n-2$, and hence $e_{n-3}(n,\ell,n-2)\leq n$. It suffices to show
$e_{n-3}(n,\ell,n-2)=n$ if $2\leq \ell \leq
\lfloor\frac{n}{2}\rfloor-1$, or $\ell=\lfloor\frac{n}{2}\rfloor$
and $n$ is odd. Let $G$ be a graph of order $n$ such that
$sdiam_{n-3}(G)\leq n-2$ and $\Delta(G)=\ell$, where $2\leq \ell
\leq \lfloor\frac{n}{2}\rfloor-1$, or
$\ell=\lfloor\frac{n}{2}\rfloor$ and $n$ is odd. Then we have the
following claim.

{\bf Claim 1.} $G$ is not a tree.

\noindent{\bf Proof of Claim 1.} Assume, to the contrary, that $G$
is a tree. Since $sdiam_{n-3}(G)\leq n-2$, it follows from Lemma
\ref{lem4-2} that $G$ contains at most two cut vertices. Then
$G=K_{1,n-1}$ or $G$ is a double star of order $n$. If
$G=K_{1,n-1}$, then $\Delta(G)=n-1>\lfloor\frac{n}{2}\rfloor$, a
contradiction. Suppose that $G$ is a double star of order $n$. Let
$u,v$ be the two centers of $G$. Then $d_G(u)\geq \lceil
\frac{n}{2}\rceil$ or $d_G(v)\geq \lceil \frac{n}{2}\rceil$, and
hence $\ell\geq \lceil \frac{n}{2}\rceil$, a contradiction. \qed

From Claim 1, $G$ is not a tree. Then $e(G)\geq n$, and hence
$e_{n-3}(n,\ell,n-2)=n$ if $2\leq \ell \leq
\lfloor\frac{n}{2}\rfloor-1$, or $\ell=\lfloor\frac{n}{2}\rfloor$
and $n$ is odd.

We now show that if $\lfloor\frac{n}{2}\rfloor+1\leq \ell \leq n-1$,
or $\ell=\lfloor\frac{n}{2}\rfloor$ and $n$ is even, then
$e_{n-3}(n,\ell,n-2)=n-1$. Let $G$ be a double star of order $n$
such that $d_G(u)=\ell$ and $d_G(v)=n-\ell$, where $u,v$ are the two
centers of $G$. From Lemma \ref{lem4-2}, we have $sdiam_{n-3}(G)\leq
n-2$. Since $\lfloor\frac{n}{2}\rfloor+1\leq \ell \leq n-1$, or
$\ell=\lfloor\frac{n}{2}\rfloor$ and $n$ is even, it follows that
$\Delta(G)=\ell$, and hence $e_{n-3}(n,\ell,n-2)\leq n-1$. So, we
have $e_{n-3}(n,\ell,n-2)=n-1$. \qed
\end{pf}

\vskip0.3cm

A graph is said to be \emph{minimally $k$-connected} if it is
$k$-connected but omitting any of the edges the resulting graph is
no longer $k$-connected.
\begin{figure}[!hbpt]
\begin{center}
\includegraphics[scale=0.8]{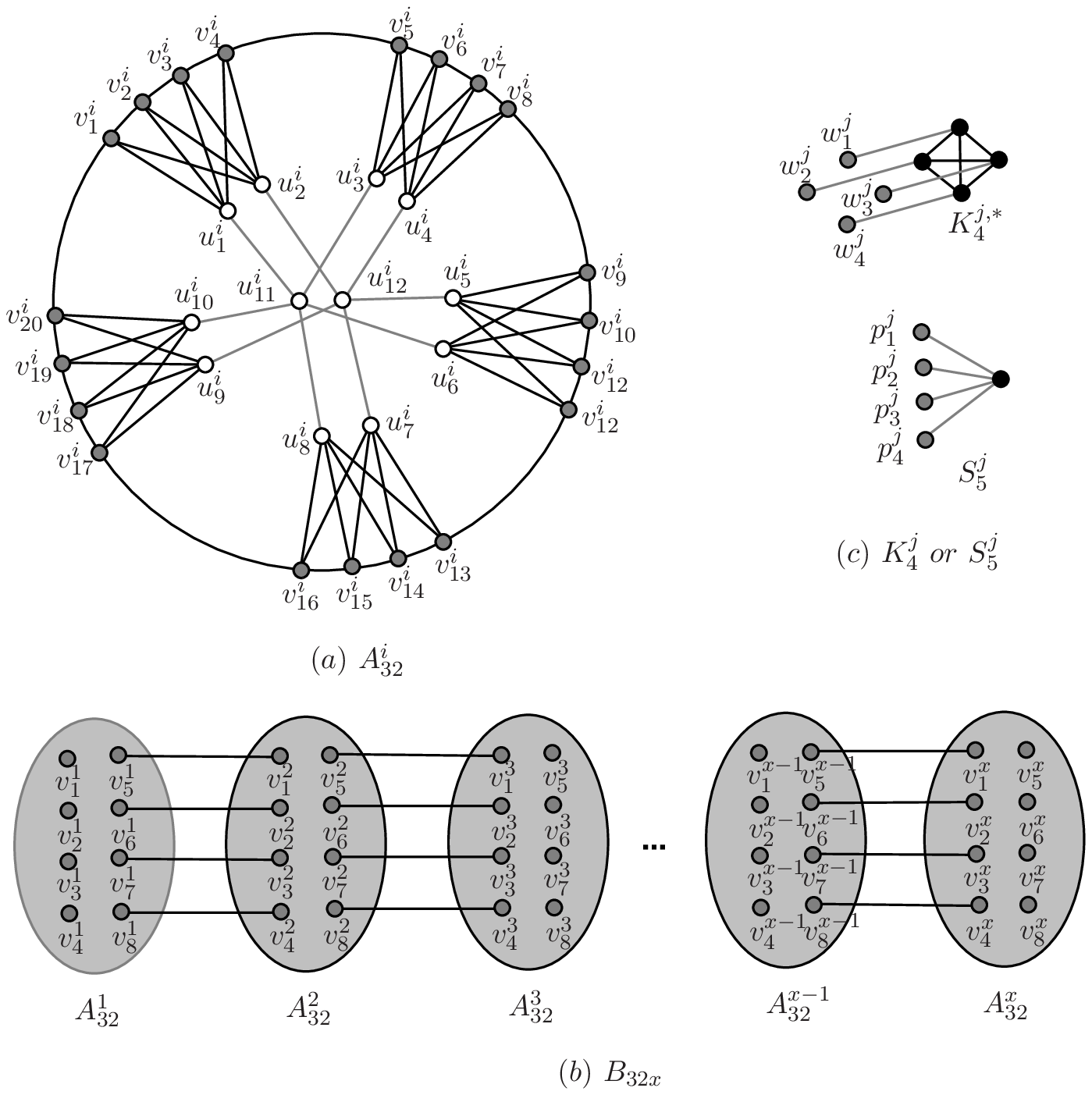}\\[0.5cm]
Figure 4: Graphs for Proposition \ref{pro4-2}.
\end{center}\label{fig7}
\end{figure}

Let $A_{32}$ be a minimally $4$-connected graph shown in Figure 4
$(a)$ (see \cite{Bollobas}, Page 18). We now give a graph $H_n$ of
order $n \ (n\geq 96)$ such that $\Delta(H_n)=\ell$ and
$sdiam_{n-3}(H_n)=n-4$ constructed by the following steps.
\begin{itemize}
\item[] \textbf{Step 1}: For each $i \ (1\leq i\leq x)$, we let $A_{32}^i$ be the copy of
$A_{32}$, where $n=32x+y$, $x=\lfloor n/32\rfloor$, and $0\leq y\leq
31$. Let $V(A_{32}^i)=\{u^i_{j}\,|\,1\leq j\leq 12\}\cup
\{v^i_{j}\,|\,1\leq j\leq 20\}$ such that $d_G(u^i_{j})=5$ for
$1\leq j\leq 12$, and $d_G(v^i_{j})=4$ for $1\leq j\leq 20$; see
Figure 4 $(a)$. Let $B_{32x}$ be a graph obtained from $A_{32}^i \
(1\leq i\leq x)$ by adding the edges in
$\{v_5^{i}v_1^{i+1}\,|\,1\leq i\leq x-1\}\cup
\{v_6^{i}v_2^{i+1}\,|\,1\leq i\leq x-1\}\cup
\{v_7^{i}v_3^{i+1}\,|\,1\leq i\leq x-1\}\cup
\{v_8^{i}v_4^{i+1}\,|\,1\leq i\leq x-1\}$; see Figure 4 $(b)$.

\item[] \textbf{Step 2}: Let $y=4z+a$, where $z=\lfloor
y/4\rfloor$, $0\leq a\leq 3$. For each $j \ (1\leq j\leq z)$, we let
$K_4^{j}$ be the complete graph of order $4$. Furthermore, let
$K_4^{j,*}$ be the graph obtained from $K_4^j$ by adding four
pendant vertices $w_1^j,w_2^j,w_3^j,w_4^j$ with four pendant edges
such that another end vertex of each pendant edge is attached on
only one vertex in $K_4^j$; see Figure 4 $(c)$. For each $j \ (1\leq
j\leq a)$, we let $S_5^{j}$ be the star of order $5$ with its leaves
$p_1^j,p_2^j,p_3^j,p_4^j$. Since $n\geq 96$, it follows that
$A_{32}^1,A_{32}^2,A_{32}^3$ all exist. Set $S_1=\{v_j^{1}\,|\,1\leq
j\leq 4\}\cup \{v_j^{1}\,|\,9\leq j\leq 20\}\subseteq V(A_{32}^1)$,
and $S_2=\{v_j^{2}\,|\,9\leq j\leq 20\}\subseteq V(A_{32}^2)$, and
$S_3=\{v_j^{3}\,|\,9\leq j\leq 20\}\subseteq V(A_{32}^3)$. Then
$|S_1\cup S_2\cup S_3|=40$. If $n\equiv 0~(mod~32)$, then
$D_n=B_{32x}$. If $n\neq 0~(mod~32)$ and $n-32x\equiv 0~(mod~4)$,
then $D_n$ is a graph obtained from $B_{32x}$ and
$K_4^{1,*},K_4^{2,*},\ldots,K_4^{z,*}$ by identifying each vertex in
$S'=\{w_j^{i}\,|\,1\leq i\leq 4, \ 1\leq j\leq z\}$ and only one
vertex in $S_1\cup S_2\cup S_3$. Since $|S'|=4z<40=|S_1\cup S_2\cup
S_3|$, for any vertex in $S'$, we can find a vertex in $S_1\cup
S_2\cup S_3$ and then identify the two vertices. If $n\neq
0~(mod~32)$ and $n-32x\neq 0~(mod~4)$, then $D_n$ is a graph
obtained from $B_{32x}$, $K_4^{1,*},K_4^{2,*},\ldots,K_4^{z,*}$ and
$S_5^{1},S_5^{2},\ldots,S_5^{a}$ by identifying each vertex in
$S'=\{w_j^{i}\,|\,1\leq i\leq 4, \ 1\leq j\leq z\}\cup
\{p_j^{i}\,|\,1\leq i\leq 4, \ 1\leq j\leq a\}$ and only one vertex
in $S_1\cup S_2\cup S_3$. Since $|S'|=4z+4a\leq 28+12=40=|S_1\cup
S_2\cup S_3|$, for any vertex in $S'$, we can find a vertex in
$S_1\cup S_2\cup S_3$ and then identify the two vertices.

\item[] \textbf{Step 3}: Let $H_n$ be the graph $D_n$ by adding
$\ell-5$ edges between $u_{12}^{1}$ and $V(G)-u_{12}^{1}$.
\end{itemize}

\vskip 0.3cm

We now in a position to give the upper and lower bounds of
$e_{n-3}(n,\ell,n-4)$.
\begin{pro}\label{pro4-2}
Let $\ell,n$ be two integers with $2\leq \ell \leq n-1$ and $n\geq
96$. Then
$$
2n-2-\lceil \ell/2\rceil\leq e_{n-3}(n,\ell,n-4)\leq 74\left\lfloor
\frac{n}{32}\right\rfloor+2i+\ell-9,
$$
where $n\equiv~i~(mod~32)$, $1\leq i\leq 31$.
\end{pro}
\begin{pf}
Let $G$ be a graph of order $n$ such that $sdiam_{n-3}(G)=n-4$ and
$\Delta(G)=\ell$, where $2\leq \ell \leq n-1$. Since
$\Delta(G)=\ell$, it follows that there exists a vertex $v$ in $G$
such that $d_G(v)=\ell$. Since $sdiam_{n-3}(G)=n-4$, it follows from
Lemma \ref{lem4-2} that $\delta(G)\geq \kappa(G)\geq 4$. For any
vertices in $V(G)-v$, its degree is at least $4$. Then $2e(G)\geq
\ell+4(n-1)$, and hence $e_{n-3}(n,\ell,n-3)\geq 2n-2-\lceil
\ell/2\rceil$.

It suffices to show $e_{n-3}(n,\ell,n-3)\leq 74\left\lfloor
\frac{n}{32}\right\rfloor+2i+\ell-9$, where $n\equiv~i~(mod~32)$,
and $1\leq i\leq 31$. Let $G=H_n$. Clearly, $\Delta(G)=\ell$. Since
$\kappa(G)\geq 4$, it follows from Lemma \ref{lem4-2} that
$sdiam_{n-3}(G)=n-4$, and hence $e_{n-3}(n,\ell,n-3)\leq e(G)=
74\left\lfloor \frac{n}{32}\right\rfloor+2i+\ell-9$, where
$n\equiv~i~(mod~32)$, $1\leq i\leq 31$.\qed
\end{pf}

\vskip 0.3cm

For the remain case $d=n-3$, we have the following.
\begin{pro}\label{pro4-3}
Let $\ell,n$ be two integers with $2\leq \ell \leq n-1$ and $n\geq
5$.

$(i)$ If $\lceil n/2\rceil+1\leq \ell\leq n-1$, then
$e_{n-3}(n,\ell,n-3)\leq 2n-\ell+1$.

$(ii)$ If $5\leq \ell\leq \lceil n/2\rceil$, then
$$
e_{n-3}(n,\ell,n-3)\leq (2\ell+3)\left\lfloor
\frac{n}{\ell+1}\right\rfloor+\ell+\left\{
\begin{array}{ll}
-8,&\mbox{{\rm if}~$n\equiv~0~(mod~\ell+1)$};\\
-5,&\mbox{{\rm if}~$n\equiv~1~(mod~\ell+1)$};\\
-2,&\mbox{{\rm if}~$n\equiv~2~(mod~\ell+1)$};\\
&\mbox{{\rm ~~~or}~$n\equiv~3~(mod~\ell+1)$};\\
2i-7,&\mbox{{\rm if}~$n\equiv~i~(mod~\ell+1)$},
\end{array}
\right.
$$
where $n=(\ell+1)x+i$ and $0\leq i\leq \ell$.

$(iii)$ If $2\leq \ell \leq n-1$, then $e_{n-3}(n,\ell,n-3)\geq
\max\{n-1+\lceil \ell/2\rceil,\lceil\frac{3n-\ell-3}{2}\rceil\}$.
\end{pro}
\begin{pf}
$(i)$ Let $G$ be a graph obtained from a wheel $W_{n-\ell+3}$ with
center $w$ and a star $K_{1,\ell-3}$ by identifying the center of
the star and one vertex of $W_{n-\ell+3}-w$. Note that
$W_{n-\ell+3}-w$ is a cycle of order $n-\ell+2$, say $C=u_1u_2\ldots
u_{n-\ell+2}u_1$. Let $u_1$ be the identifying vertex. Clearly,
$d_G(w)=n-\ell+2$, $d_G(u_1)=\ell$ and $d_G(u_i)=3$ for each $i \
(2\leq i\leq n-\ell+2)$. Since $\lceil n/2\rceil+1\leq \ell\leq
n-1$, it follows that $\ell\geq n-\ell+2$, and hence
$\Delta(G)=\ell$. Since $G$ contains only one cut vertex, one can
easily check that $sdiam_{n-3}(G)\leq n-3$. So, we have
$e_{n-3}(n,\ell,n-3)\leq e(G)=2n-\ell+1$.

$(ii)$ Let $n=(\ell+1)x+y$, where $x=\lfloor n/(\ell+1)\rfloor$,
$0\leq y\leq \ell$. For each $i \ (1\leq i\leq x)$, we let
$W_{\ell+1}^i$ be the wheel of order $\ell+1$, with its center
$w_i$. Note that $W_{\ell+1}^i-w_i$ is a cycle of order $\ell$, say
$v^i_1v^i_2\ldots v^i_{\ell}v^i_1$. Let $F_{(\ell+1)x}$ be a graph
obtained from $W_{\ell+1}^i \ (1\leq i\leq x)$ by adding the edges
in $\{v_j^{i}v_j^{i+1}\,|\,1\leq i\leq x-1, \ 1\leq j\leq 3\}$. We
now give a graph $D_n$ of order $n \ (n\geq \ell+1)$ constructed in
the following way.
\begin{itemize}
\item If $n\equiv~0~(mod~\ell+1)$, then $D_n=F_{(\ell+1)x}$.

\item If $n\equiv~1~(mod~\ell+1)$, then $D_n$ is the graph obtained from $F_{(\ell+1)x}$
by adding a new vertex $u$ and three edges $uv^1_1,uv^1_2,uv^1_3$.

\item If $n\equiv~2~(mod~\ell+1)$, then $D_n$ is the graph obtained from $F_{(\ell+1)x}$
by adding two new vertices $u_1,u_2$ and six edges
$\{u_1u_2,u_1v^1_1,u_1v^1_2,u_2v^1_3,u_2v^1_4,u_2v^1_5\}$.

\item If $n\equiv~3~(mod~\ell+1)$, then $D_n$ is the graph obtained from $F_{(\ell+1)x}$
by adding two new vertices $u_1,u_2,u_3$ and six edges
$\{u_1u_2,u_1u_3,u_2u_3,u_1v^1_1,u_2v^1_2,u_3v^1_3\}$.

\item If $n\equiv~y~(mod~\ell+1) \ (4\leq y\leq \ell)$, then $D_n$ is the graph obtained from
$F_{(\ell+1)x}$ and a new wheel $W^*_{y}$ with center $w^*$ and $y$
vertices by adding three edges
$\{v^*_2v^1_1,v^*_2v^1_2,v^*_{3}v^1_3\}$, where
$V(W^*_{y})-w^*=\{v^*_2,v^*_2,\ldots,v^*_{y-1}\}$.
\end{itemize}
Let $H_n$ be the graph $D_n$ by adding $\ell-5$ edges between
$v_{1}^{1}$ and $V(G)-v_{1}^{1}$. Clearly, $\Delta(H_n)=\ell$. Since
$H_n$ is $3$-connected, one can easily check that
$sdiam_{n-3}(G)\leq n-3$. Since
$$
e(G)=2\ell x+3(x-1)+\ell-5+\left\{
\begin{array}{ll}
0,&\mbox{{\rm if}~$y=0$};\\
3,&\mbox{{\rm if}~$y=1$};\\
6,&\mbox{{\rm if}~$y=2$};\\
&\mbox{{\rm ~~~or}~$y=3$};\\
2(y-1)+3,&\mbox{{\rm if}~$4\leq y\leq \ell$},
\end{array}
\right.
$$
it follows that
$$
e_{n-3}(n,\ell,n-3)\leq (2\ell+3)\left\lfloor
\frac{n}{\ell+1}\right\rfloor+\ell+\left\{
\begin{array}{ll}
-8,&\mbox{{\rm if}~$n\equiv~0~(mod~\ell+1)$};\\
-5,&\mbox{{\rm if}~$n\equiv~1~(mod~\ell+1)$};\\
-2,&\mbox{{\rm if}~$n\equiv~2~(mod~\ell+1)$};\\
&\mbox{{\rm ~~~or}~$n\equiv~3~(mod~\ell+1)$};\\
2i-7,&\mbox{{\rm if}~$n\equiv~i~(mod~\ell+1)$},
\end{array}
\right.
$$
where $n=(\ell+1)x+i$ and $0\leq i\leq \ell$.

$(iii)$ Let $G$ be a graph of order $n$ such that
$sdiam_{n-3}(G)\leq n-3$ and $\Delta(G)=\ell$, where $2\leq \ell
\leq n-1$. Then $G$ is $2$-connected, or $G$ contains only one cut
vertex. If $G$ is $2$-connected, then there exists a vertex $v$ in
$G$ such that $d_G(v)=\ell$, since $\Delta(G)=\ell$. For any
vertices in $V(G)-v$, its degree is at least $2$. Then $2e(G)\geq
\ell+2(n-1)$, and hence $e_{n-3}(n,\ell,n-3)\geq n-1+\lceil
\ell/2\rceil$. Suppose that $G$ contains only one cut vertex, say
$v$. Then each connected component of $G\setminus v$ is a connected
subgraph of order at least $3$, or an edge of $G$, or an isolated
vertex. Let $w_1,w_2,\ldots,w_r$ be the isolated vertices,
$e_1,e_2,\ldots,e_s$ be the edges, and $C_1,C_2,\ldots,C_t$ be the
connected components of order at least $3$ in $G\setminus v$. Then
we have the following claim.

{\bf Claim 1.} For any $w\in \bigcup_{i=1}^tV(C_i)$, if $d_G(w)=2$,
then $wv\in E(G)$.

\noindent{\bf Proof of Claim 1.} Assume, to the contrary, that
$wv\notin E(G)$ for any $w\in \bigcup_{i=1}^tV(C_i)$ with
$d_G(w)=2$. Without loss of generality, let $w\in V(C_1)$. Then
there exist two vertices $u_1,u_2$ in $C_1$ such that $u_1w\in
E(C_1)$ and $u_2w\in E(C_1)$. Choose $S\subseteq V(G)$ with
$|S|=n-3$ such that $w\in S$ but $u_1,u_2,v\notin S$. Then any
$S$-Steiner tree must occupy $v$ and one of $u_1,u_2$, and hence
$d_G(S)\geq n-2$, a contradiction. \qed

\vskip 0.3cm

From Claim 1, we suppose that there are $x$ vertices
$u_1,u_2,\ldots,u_x$ in $\bigcup_{i=1}^tV(C_i)$ such that its degree
is $2$. Then for any vertex in $(\bigcup_{i=1}^tV(C_i))\setminus
\{u_1,u_2,\ldots,u_x\}$, its degree is at least $3$.

If $d_G(v)=\ell$, then $x\leq \ell-r-2s$ and $r\leq \ell$.
Furthermore, we have
\begin{eqnarray*}
2e(G)&\geq &r+4s+\ell+2x+3(n-r-2s-x-1)\\
&=&3n-2r-2s-x+\ell-3\\
&\geq&3n-2r-2s-(\ell-r-2s)+\ell-3=3n-r-3=3n-\ell-3,
\end{eqnarray*}
and hence $e(G)\geq \lceil\frac{3n-\ell-3}{2}\rceil$.

If $d_G(v)\neq \ell$, then there exists a vertex $u\in V(G)-v$ such
that $d_G(u)=\ell$. Clearly, $r+2s+t\leq d_G(v)\leq \ell$. Since
$x+2s+r\leq \ell-1$, it follows that
\begin{eqnarray*}
2e(G)&\geq &r+4s+d_G(v)+\ell+2x+3(n-r-2s-x-2)\\
&\geq&r+4s+(r+2s+t)+\ell+2x+3(n-r-2s-x-2)\\
&=&3n-r-x+t+\ell-6\\
&\geq&3n-(\ell-1-2s)+t+\ell-6\\
&=&3n+2s+t-5\geq 3n-5,
\end{eqnarray*}
and hence $e(G)\geq \lceil\frac{3n-5}{2}\rceil$.

From the above argument, we conclude that $e_{n-3}(n,\ell,n-3)\geq
\max\{n-1+\lceil \ell/2\rceil,\lceil\frac{3n-\ell-3}{2}\rceil\}$.
\qed
\end{pf}

From Propositions \ref{pro4-1}, \ref{pro4-2} and \ref{pro4-3}, we
have the following theorem.
\begin{thm}\label{th4-1}
$(i)$ For $2\leq \ell \leq n-1$, $e_{n-1}(n,\ell,n-1)=n-1$.

$(ii)$ For $2\leq \ell \leq n-1$ and $n\geq 4$,
$$
e_{n-3}(n,\ell,n-2)=\left\{
\begin{array}{ll}
n,&\mbox{{\rm if}~$2\leq \ell \leq \lfloor\frac{n}{2}\rfloor-1$};\\
&\mbox{{\rm ~~~or}~$\ell=\lfloor\frac{n}{2}\rfloor$~and~n~is~odd};\\
n-1,&\mbox{{\rm if}~$\lfloor\frac{n}{2}\rfloor+1\leq \ell \leq n-1$};\\
&\mbox{{\rm ~~~or}~$\ell=\lfloor\frac{n}{2}\rfloor$~and~n~is~even}.\\
\end{array}
\right.
$$

$(iii)$ For $2\leq \ell \leq n-1$ and $n\geq 96$,
$$
2n-2-\lceil \ell/2\rceil\leq e_{n-3}(n,\ell,n-4)\leq 74\left\lfloor
\frac{n}{32}\right\rfloor+2i+\ell-9,
$$
where $n\equiv~i~(mod~32)$, $1\leq i\leq 31$.

$(iv)$ If $2\leq \ell \leq n-1$, then $e_{n-3}(n,\ell,n-3)\geq
\max\{n-1+\lceil \ell/2\rceil,\frac{3n-\ell-5}{2}\}$. If $\lceil
n/2\rceil+1\leq \ell\leq n-1$, then $e_{n-3}(n,\ell,n-3)\leq
2n-\ell+1$. If $5\leq \ell\leq \lceil n/2\rceil$, then
$$
e_{n-3}(n,\ell,n-3)\leq (2\ell+3)\left\lfloor
\frac{n}{\ell+1}\right\rfloor+\ell+\left\{
\begin{array}{ll}
-8,&\mbox{{\rm if}~$n\equiv~0~(mod~\ell+1)$};\\
-5,&\mbox{{\rm if}~$n\equiv~1~(mod~\ell+1)$};\\
-2,&\mbox{{\rm if}~$n\equiv~2~(mod~\ell+1)$};\\
&\mbox{{\rm ~~~or}~$n\equiv~3~(mod~\ell+1)$};\\
2i-7,&\mbox{{\rm if}~$n\equiv~i~(mod~\ell+1)$},
\end{array}
\right.
$$
where $n=(\ell+1)x+i$ and $0\leq i\leq \ell$.
\end{thm}

\section{For general $k$}

In \cite{MaoMelekianCheng}, Mao et al. obtained the following
result.
\begin{lem}\label{lem5-1}
Let $\ell,n$ be two integers with $1\leq \ell\leq n-2$, and let $G$
be a graph of order $n$. Then $\kappa(G)\geq \ell$ if and only if
$sdiam_{n-\ell+1}(G)=n-\ell$.
\end{lem}

In this section, we construct a graph and give an upper bound of $e_k(n,\ell,d)$
for general $k$, $\ell$, and $d$.
\begin{thm}\label{th5-1}
Let $k,\ell,d$ be three integers with $2\leq k\leq n$, $2\leq
\ell\leq n-1$, and $k-1\leq d\leq n-1$.

$(i)$ If $d=k-1$, $\lceil \frac{n+1}{2}\rceil\leq k\leq n$, and
$\max\{n-k+1,\lceil\frac{n}{2}\rceil\}<\ell\leq n-1$, then
$$
\left\lceil \frac{\ell+(n-1)(n-k+1)}{2}\right\rceil \leq
 e_k(n,\ell,d)\leq \frac{(n-1)^2}{4}+\ell.
$$

$(ii)$ If $2\leq k\leq d$, $k\leq d\leq n-1$, and $2+\lceil
\frac{n-d+k-3}{d-k+1}\rceil\leq \ell\leq n-1$, then
$$
e_k(n,\ell,d)=n-1.
$$
\end{thm}
\begin{pf}
$(i)$ We first consider the lower bound. From Lemma \ref{lem5-1},
for a connected graph $G$ of order $n$, $sdiam_{k}(G)=k-1$ if and
only if $\kappa(G)\geq n-k+1$. Let $G$ be a graph of order $n$ such
that $sdiam_{k}(G)=k-1$ and $\Delta(G)=\ell$, where $2\leq \ell \leq
n-1$. Since $sdiam_{k}(G)=k-1$, it follows that $\delta(G)\geq
\kappa(G)\geq n-k+1$. Since $\Delta(G)=\ell$, it follows that there
exists a vertex $v$ in $G$ such that $d_G(v)=\ell$. For any vertex
in $V(G)-v$, its degree is at least $n-k+1$. Then $2e(G)\geq
\ell+(n-1)(n-k+1)$, and hence $e_{k}(n,\ell,k-1)\geq \lceil
\frac{\ell+(n-1)(n-k+1)}{2}\rceil$.

Next, we consider the upper bound. For
$\max\{n-k+1,\lceil\frac{n}{2}\rceil\}\leq \ell\leq n-1$, we let
$K_{a,b}$ be a complete bipartite graph of order $n=a+b$ with $a\geq
b$. Let $U,V=\{v_1,v_2,\ldots,v_b\}$ be the parts of order $a,b$ in
$K_{a,b}$. Let $G$ be a graph obtained from $K_{a,b}$ by adding
edges $v_1v_i \ (2\leq i\leq \ell-a)$, where $\ell\leq a+b-1=n-1$.
Then $\Delta(G)=a+(\ell-a)=\ell$. Since
$\max\{n-k+1,\lceil\frac{n}{2}\rceil\}<\ell\leq n-1$ and $a\geq b$,
it follows that for any $S\subseteq V(G)$ and $|S|=k$, we have
$S\cap U\neq \emptyset$ and $S\cap U\neq \emptyset$, and hence
$sdiam_{k}(G)=k-1$. So $e_k(n,\ell,d)\leq
ab+\ell-a=a(n-a)-a+\ell\leq \frac{(n-1)^2}{4}+\ell$.

$(ii)$ Let $P_{d-k+3}=v_1v_2\ldots v_{d-k+3}$ be a path of order
$d-k+3$. Set $x=\lceil \frac{n-d+k-3}{d-k+1}\rceil$,
$U_i=\{u_{i,1},u_{i,2},\ldots,u_{i,x}\,|\,2\leq i\leq d-k+1\}$, and
$U_{d-k+2}=\{u_{d-k+2,1},u_{d-k+2,2},\ldots,u_{d-k+2,p}\}$, where
$p=n-(d-k)(x+1)+3$. Let $T'$ be a tree of order $n$ obtained from
$P_{d-k+3}$ by adding all the vertices in
$\bigcup_{i=2}^{d-k+2}U_i$, and then adding the edges in
$\bigcup_{i=2}^{d-k+2}E_{T'}[v_i,U_{i}]$. Let $T$ be a tree obtained
from $T'$ by gradually deleting $\ell-x-2$ vertices in
$V(T')\setminus (\{v_1,v_2\}\bigcup U_2)$ and then adding $\ell-x-2$
pendant edges at $v_2$. Clearly, $d_T(v_2)=\ell=\Delta(T)$. Since
$sdiam_{k}(T)\leq d$, it follows that $e_k(n,\ell,d)=n-1$.
\end{pf}


\begin{thebibliography}{1}

\bibitem{Ali}
P. Ali, The Steiner diameter of a graph with prescribed girth, {\it
Discrete Math.} 313(12)(2013), 1322--1326.


\bibitem{AliDM}
P. Ali, P. Dankelmann, and S. Mukwembi, Upper bounds on the Steiner
diameter of a graph, {\it Discrete Appl. Math.} 160(12) (2012),
1845--1850.


\bibitem{Bollobas}
B. Bollob\'{a}s, {\it Extremal Graph Theory}, Acdemic press, 1978.

\bibitem{B26}
B. Bollob\'{a}s, Graphs with a given diameter and maximal valency
and with a minimal number of edges, in:``Combinatorial Mathematics
and its Applications'' (Welsh, D.J.A., ed.) {\it Academic Press}
London and New York, 1971), 25--37.

\bibitem{Bondy}
J.A. Bondy and U.S.R. Murty, Graph Theory, GTM 244, Springer, 2008.

\bibitem{Harary}
F. Buckley and F. Harary, Distance in Graphs, Addision-Wesley,
Redwood City, CA, 1990.

\bibitem{Caceresa}
J. C\'{a}ceresa, A. M\'{a}rquezb, and M. L. Puertasa, Steiner
distance and convexity in graphs, {\it European J. Combin.}
29(2008)(3), 726--736.

\bibitem{Chartrand}
G. Chartrand, O.R. Oellermann, S. Tian, and H.B. Zou, Steiner
distance in graphs, {\it \'{C}asopis pro p\v{e}stov\'{a}n\'{i}
matematiky} 114(1989), 399--410.

\bibitem{ChartrandOZ}
G. Chartrand, F. Okamoto, and P. Zhang, Rainbow trees in graphs and
generalized connectivity, {\it Networks} 55(4)(2010), 360--367.

\bibitem{Chung}
F.R.K. Chung, Diameter of graphs: old problems and new results, 18th
Southeastern Conf. on Combinatorics, Graph Theory and Computing,
1987.

\bibitem{DankelmannOS} P. Dankelmann, O.R. Oellermann, and H.C. Swart,
The average Steiner distance of a graph, {\it J. Graph Theory}
22(1)(1996), 15--22.

\bibitem{DankelmannSO}
P. Dankelmann, H.C. Swart, and O.R. Oellermann, On the average
Steiner distance of graphs with prescribed properties, {\it Discrete
Appl. Math.} 79(1-3)(2008), 91--103.

\bibitem{DankelmannSO2}
P. Dankelmann, H. Swart, and O.R. Oellermann, Bounds on the Steiner
diameter of a graph, {\it Combinatorics, Graph Theory, and
Algorithms, Vol. I, II, New Issues Press}, Kalamazoo, MI, 1999.

\bibitem{DayOS}
D.P. Day, O.R. Oellermann, and H.C. Swart, Steiner
distance-hereditary Graphs, {\it SIAM J. Discrete Math.} 7(3)
(1994), 437--442.

\bibitem{Du}
D.Z. Du, Y.D. Lyuu, and D.F. Hsu, Line digraph iteration and
connectivity analysis of de Bruijn and Kautz graphs, {\it IEEE
Trans. Comput.} 42(5)(1994), 612--616.

\bibitem{Dobrynin} A. Dobrynin, R. Entringer, and I. Gutman,
Wiener index of trees: theory and application, {\it Acta Appl.
Math.} 66(2001), 211--249.

\bibitem{Moscarini}
A. D'Atri and M. Moscarini, Distance-Hereditary Graphs, Steiner
Trees, and Connected Domination, {\it SIAM J. Comput.} 17(3)(1988),
521--538.

\bibitem{ER4}
P. Erd\"{o}s and A. R\'{e}nyi, On a problem in the theory of graphs
(in Hungarian), {\it Publ. Math. Inst. Hungar. Acad. Sci.} 7(1962).

\bibitem{ERS1}
P. Erd\"{o}s, A. R\'{e}nyi, and V.T. S\'{o}s, On a problem of graph
theory, {\it Studia Sci. Math. Hungar.} 1(1966), 215--235.

\bibitem{FurtulaGK} B. Furtula, I. Gutman, and V. Katani\'{c},
Three-center Harary index and its applications, {\it Iranian J.
Math. Chem.} 7(1)(2016), 61--68.

\bibitem{Meyer}
F.J. Meyer, D.K. Pradhan, Flip trees, {\it IEEE Trans. Computers}
37(3)(1987), 472--478.

\bibitem{GareyJ}
M.R. Garey and D.S. Johnson, Computers and Intractibility: A Guide
to the Theory of NP-Completeness, {\it IEEE Trans. Computers},
Freeman \& Company, New York, 1979.

\bibitem{GoddardOS}
W. Goddard, O.R. Oellrmann, and H.C. Swart, Steiner distance stable
graphs, {\it Discrete Math.} 132(1-3)(1994), 65--73.

\bibitem{Goddard}
W. Goddard and O.R. Oellrmann, Distance in graphs, in: M. Dehmer
(Ed.) Structural Analysis of Complex Networks, Birkh\"{a}user,
Dordrecht, 2011, 49--72.

\bibitem{GutmanSDD} I. Gutman,
On Steiner degree distance of trees,
{\it Appl. Math. Comput.\/} 283(2016), 163--167.

\bibitem{GFL} I. Gutman, B. Furtula, and X. Li,
Multicenter Wiener indices and their applications, {\it J. Serb.
Chem. Soc.} 80 (2015) 1009--1017.

\bibitem{Hakimi}
S.L. Hakimi, Steiner's problem in graph and its implications {\it
Networks}, 1(2)(1971), 113--133.

\bibitem{HO1}
N.P. Homenko and N.A. Ostroverhii, Diameter-critical graphs (in
Russian), {\it Ukrainian Math. J.} 22(1970), 637--646.

\bibitem{HS5}
N.P. Homenko and V.V. Strok, Some combinatorial indentities for sums
of composition coefficient (in Russian), {\it Ukrainian Math. J.} 23
(1971), 830--837.

\bibitem{HwangRW}
F.K. Hwang, D.S. Richards, and P. Winter, The Steiner Tree Problem,
North-Holland, Amsterdam, 1992.

\bibitem{Hsu}
D.F. Hsu, On container width and length in graphs, groups, and
networks, {\it IEICE Transaction on Fundamentals of Electronics,
Communications and Computer Science}, E77-A (1994), 668--680.

\bibitem{Hsu2}
D.F. Hsu and T. {\L}uczak, Note on the $k$-diameter of $k$-regular
$k$-connected graphs, {\it Discrete Math.} 133(1-3)(1994), 291--296.

\bibitem{Levi}
A.Y. Levi, Algorithm for shortest connection of a group of graph
vertices, {\it Sov. Math. Dokl.} 12(1971), 1477--1481.

\bibitem{LMG} X. Li, Y. Mao, and I. Gutman,
The Steiner Wiener index of a graph, {\it Discuss. Math. Graph
Theory} 36(2016), 455--465.

\bibitem{LMG2} X. Li, Y. Mao, and I. Gutman,
      Inverse problem on the Steiner Wiener index,
        {\it Discuss. Math. Graph Theory}, in press.

\bibitem{Mao}
Y. Mao, The Steiner diameter of a graph, {\it Bull. Iran. Math.
Soc.}, in press.

\bibitem{Mao2}
Y. Mao, Steiner $3$-diameter, maximum degree and size of a graph,
arXiv:1703.04974 [math.CO] 2017.

\bibitem{MaoDas}
Y. Mao and K.C. Das, Steiner Gutman index, {\it MATCH Commun. Math.
Comput. Chem.}, in press.

\bibitem{MaoMelekianCheng}
Y. Mao, C. Melekian, and E. Cheng, The Steiner $(n-3)$-diameter of a
graph, arXiv:1703.03984 [math.CO] 2017.


\bibitem{MWG}
Y. Mao, Z. Wang, and I. Gutman,
      Steiner Wiener index of graph products,
      {\it Trans. Combin.} 5(3)(2016), 39--50.

\bibitem{MWGK}
Y. Mao, Z. Wang, I. Gutman, and A. Klobu\v{c}ar,
      Steiner degree distance,
      {\it MATCH Commun. Math. Comput. Chem.} 78(1)(2017), 221--230.

\bibitem{MWGL}
Y. Mao, Z. Wang, I. Gutman, and H. Li, Nordhaus-Gaddum-type results
for the Steiner Wiener index of graphs, {\it Discrete Appl. Math.}
219(2017) 167--175.



\bibitem{Meyer}
F.J. Meyer and D.K. Pradhan, Flip trees, {\it IEEE Trans. Computers}
37(1987)(3), 472--478.

\bibitem{Oellermann} O.R. Oellermann, From Steiner centers to
Steiner medians, {\it J. Graph Theory} 20(2)(1995), 113--122.

\bibitem{Oellermann2} O.R. Oellermann,
On Steiner centers and Steiner
medians of graphs, {\it Networks} 34(1999), 258--263.

\bibitem{OellermannT}
O.R. Oellermann and S. Tian, Steiner Centers in graphs, {\it J.
Graph Theory} 14(5)(1990), 585--597.

\bibitem{O6}
O. Ore, Diameter in graphs, {\it J. Combin. Theory} \textbf{5}
(1968), 75--81.

\bibitem{Rouv1} D.H. Rouvray,
Harry in the limelight: The life and times of Harry Wiener, in: D.
H. Rouvray, R. B. King (Eds.), \emph{Topology in Chemistry --
Discrete Mathematics of Molecules}, Horwood, Chichester, 2002, pp.
1--15.

\bibitem{Rouv2} D.H. Rouvray,
The rich legacy of half century of the Wiener index, in: D.H.
Rouvray, R.B. King (Eds.), \emph{Topology in Chemistry -- Discrete
Mathematics of Molecules}, Horwood, Chichester, 2002, pp. 16--37.


\bibitem{WangMLY}
Z. Wang, Y. Mao, H. Li, and C. Ye, The Steiner $4$-diameter of a
graph, arXiv:1702.05681 [math.CO] 2017.

\bibitem{W11}
M.E. Watkins, A lower bound for the number of vertices of a graph,
{\it Amer. Math. Monthly} \textbf{74} (1976), 297.

\bibitem{Xu}
K. Xu, M. Liu, K.C. Das, I. Gutman, and B. Furtula, A survey on
graphs extremal with respect to distance--based topological indices,
{\it MATCH Commun. Math. Comput. Chem.} 71(2014), 461--508.
\end{thebibliography}
\end{document}